\documentclass{article}

\usepackage{authblk}
\usepackage{cite}
\usepackage{epsf,latexsym,amssymb}
\usepackage{graphicx}
\usepackage{amsmath}  
\usepackage{color}

\newcommand{\R}{{\mathbb R}}  

\newcommand{\be}[1]{\begin{equation}\label{#1}}
\newcommand{\ee}{\end{equation}}
\newcommand{\halmos}{\rule{1ex}{1.4ex}}
\newcommand{\mybox}{\hfill $\Box$} 
\newcommand{\beqn}{\begin{eqnarray*}}
\newcommand{\eeqn}{\end{eqnarray*}}
\newcommand{\beqnum}{\begin{eqnarray}}
\newcommand{\eeqnum}{\end{eqnarray}}
\newcommand{\ben}{\begin{enumerate}}
\newcommand{\een}{\end{enumerate}}
\newcommand{\mypmatrix}[1]{\left(\begin{array}{cccccccccccc}#1\end{array}\right)}

\newcommand{\x}{\xi } 

\newtheorem{theorem}{Theorem}
\newtheorem{itlemma}{Lemma}[section] 
\newtheorem{itproposition}[itlemma]{Proposition}
\newtheorem{itcorollary}[itlemma]{Corollary}
\newtheorem{itremark}[itlemma]{Remark}
\newtheorem{itdefinition}[itlemma]{Definition}

\newenvironment{lemma}{\begin{itlemma}\rm}{\end{itlemma}} 
\newenvironment{remark}{\begin{itremark}\rm}{\end{itremark}} 
\newenvironment{corollary}{\begin{itcorollary}\rm}{\end{itcorollary}}
\newenvironment{proposition}{\begin{itproposition}\rm}{\end{itproposition}}
\newenvironment{definition}{\begin{itdefinition}\rm}{\end{itdefinition}}

\newenvironment{proof}{\noindent {\em Proof}.\ }{\hspace*{\fill}$\halmos$\medskip}

\newcommand{\bl}[1]{\begin{lemma}\label{#1}}
\newcommand{\br}[1]{\begin{remark}\label{#1}}
\newcommand{\bt}[1]{\begin{theorem}\label{#1}}
\newcommand{\bd}[1]{\begin{definition}\label{#1}}
\newcommand{\bp}[1]{\begin{proposition}\label{#1}}
\newcommand{\bc}[1]{\begin{corollary}\label{#1}}
\newcommand{\ec}{\mybox\end{corollary}}
\newcommand{\el}{\end{lemma}}
\newcommand{\er}{\mybox\end{remark}}
\newcommand{\et}{\end{theorem}}
\newcommand{\ed}{\mybox\end{definition}}
\newcommand{\ep}{\mybox\end{proposition}}
\newcommand{\epr}{\end{proof}}
\newcommand{\bpr}{\begin{proof}}

\newcommand{\diag}{\mbox{diag}\,}

\newcommand{\ab}{r}
\newcommand{\bb}{s}
\newcommand{\st}{\, | \,}

\DeclareMathOperator{\spn}{span}

\def\twodigits#1{\ifnum#1<10 0\fi\the#1}
\makeatletter
  \newcommand*\short[1]{\expandafter\@gobbletwo\number\numexpr#1\relax}
\makeatother
\usepackage{datetime}

\title{An observability result related to active sensing }
\author[1]{Eduardo D. Sontag}
\author[2]{Debojyoti Biswas}
\author[2,3]{Noah J. Cowan}
\affil[1]{Department of Electrical and Computer Engineering and Department of Bioengineering, Northeastern University, Boston, Massachusetts, United States}
\affil[2]{Laboratory for Computational Sensing and Robotics, Johns Hopkins University, Baltimore, Maryland, United States.}
\affil[3]{Department of Mechanical Engineering, Johns Hopkins University,\newline Baltimore, Maryland, United States.}

\begin{document}
\maketitle



\medskip
\section{Abstract}
For a general class of translationally invariant systems with a specific category of nonlinearity in the output, this paper presents necessary and sufficient conditions for global observability. Critically, this class of systems cannot be stabilized to an isolated equilibrium point by dynamic output feedback. These analyses may help explain the active sensing movements made by animals when they perform certain motor behaviors, despite the fact that these active sensing movements appear to run counter to the primary motor goals. The findings presented here establish that active sensing underlies  the maintenance of observability for such biological systems, which are inherently nonlinear due to the presence of the high-pass sensor dynamics.

\section{Introduction}
Active sensing is the process of expending energy, typically through
movement, for the purpose of sensing
\cite{gibson1962observations,bajcsyactive1988,schroeder2010dynamics}. Animals
use this strategy to enhance sensory information across sensory
modalities e.g.,\ echolocation
\cite{nelson2006sensory,wohlgemuth2016action}, whisking
\cite{mitchinson2011active,bush2016whisking} and other forms of touch
\cite{prescott2011active,saig2012motor}, electrosense
\cite{hofmann2013sensory,chen2020tuning,stamper2012active}, and vision
\cite{ahissar2012seeing,michaiel2020dynamics}. It is well established
that conditions of decreased sensory acuity leads to increased active
movements \cite{lockey2015one, sponberg2015luminance,
  stockl2017comparative, deora2021tactile, stamper2012active,
  chen2020tuning, catania2013stereo, michaiel2020dynamics,
  wohlgemuth2016action, kiemel2002multisensory, rucci2015unsteady} but
its actual role in relation to task-level control remains
underexplored.  The ubiquity of active sensing in nature motivates us
to explore the mathematical conditions that might necessitate active
sensing. Our theory is that active sensing is at least in part borne
out of the needs of nonlinear state estimation. We hypothesize that
animals---through active sensing---generate time-varying motor
commands that continuously stimulate their sensory receptors so that
the system states can be estimated with satisfactory error bounds from
the sensor measurements. In essence, these movements aim to maintain
the observability of the system.

A dominant paradigm in control systems engineering involves designing
state feedback and state estimation independently, an approach can be
applied successfully to a wide range of system designs.  Indeed, for
linear plants corrupted by Gaussian noise, there is a separation
principle: it is not only satisfactory to separate state estimation
from the task-level control design, but, in fact, it is \emph{optimal}
to perform this decomposition. In particular, the
linear-quadratic-Gaussian (LQG) controller decomposes into a
linear-quadratic regulator (LQR) applied to the optimal state estimate
which comes from a Kalman filter. Critically, the Kalman filter does
not depend on the LQR cost function, and the LQR gains do not depend
on the sensor noise and process noise.  Conceptually speaking,
``active sensing'' is the opposite approach to applying a separation
principle: control inputs are specifically designed to excite sensory
receptors, presumably in service to the state estimator. This may be,
at least in part, because biological sensory systems often stop
responding to persistent (i.e. ``DC'') stimuli, via sensory
``adaptation'' \cite{22_annual_reviews_tutorial_imp,2019_deepak_integral_controller,imp03,andrews2008approximate,shoval2011symmetry} or ``perceptual fading''
\cite{hofmann2014motor,jun2016active}.

In this paper, we formalize a class of nonlinear systems that have a
simple high-pass sensory output that mimics sensory adaptation or
perceptual fading. Under some simplified modeling assumptions,
reviewed below, this implies that linear observability is lost, which
means the usual LQG-style framework does not apply. However, under
some interesting modeling conditions, nonlinear observability
persists. Critically, nonlinear observability does not necessarily
afford a separation principle: the control signal may need to contain
ancillary energy that is expressly for the purpose of state
estimation, and may be in conflict with task goals. Indeed, the energy
expended for active sensing movements do not necessarily directly
serve a motor control goal, and are instead believed to improve
sensory feedback and prevent perceptual fading
\cite{stamper2012active,hofmann2014motor,jun2016active}.

The organization of the paper is as follows. Section \ref{sec:bio}
motivates the model structure from prior work and Section
\ref{sec:mainresult} generalizes the model and presents the main
theorem. Section \ref{sec:proof} has the proof of the main
theorem. The Appendix provides some background concepts for the ease
of understanding the tools used in Section \ref{sec:bio} and \ref{sec:proof}.


\section{Biological motivation and simplified system} \label{sec:bio}

Station keeping behavior in weakly electric fish, \textit{Eigenmannia
  virescens}, provides an ideal system for investigating the interplay
between active sensing and task-level control
\cite{stamper2012active,stamperusing2019,biswas2018closed,chen2020tuning}. These
fish routinely maintain their position relative to a moving refuge
and uses both vision and electrosense to collect the necessary sensory
information from its environment \cite{cowan2007critical,
  cowan2014feedback,sutton2016dynamic,roselongitudinal1993}. While
tracking the refuge position (i.e., task-level control), the fish
additionally produce rapid ``whisking-like'' forward and backward
swimming movements (i.e., active sensing). When vision is limited (for
example, in darkness), the fish increase their active sensing
movements \cite{stamper2012active,biswas2018closed}. This increased
motion compensates the lack of visual cues \cite{chen2020tuning}.

Suppose $x$ is the position of an animal and $z = \dot x$ is its
velocity as it moves in one degree of freedom. We assume that a
sensory receptor measures only the local rate of change of a stimulus,
$s(x)$ as the animal moves relative to the sensory scene,
i.e. $y=\dfrac{d}{dt}s(x)$. Defining $\gamma(x) := \dfrac{d}{dx}s(x)$, we arrive
at a $2$-dimensional, single-input, single-output normalized
mass-damper system of the following form
\cite{sefati2013mutually,kunapareddy2018recovering}:
\begin{equation}
\label{eq:simple-system}
\begin{aligned}
\dot x &= z, \quad && x\in \R\\
\dot z &= -z \,+\, u, \quad && z,u\in \R\\
y   &= \dfrac{d}{dt}s(x)=\gamma(x)\,z, \quad && y\in \R
\end{aligned}
\end{equation}
where the mass and the damping constant both are assumed to be unity. Linearization of the above system (\ref{eq:simple-system})  around any equilibrium, $(x^*, 0)$, is given by $(A, B, C)$ as follows:
\begin{align*}
A = \begin{bmatrix}
0 &1\\
0 &-1
\end{bmatrix}, \quad B = \begin{bmatrix}
0 \\
1
\end{bmatrix}, \quad C = \begin{bmatrix}
0 & \gamma^*
\end{bmatrix},
\end{align*}
where $\gamma^* = \gamma(x^*)$. Clearly (A,C) is not observable irrespective
of $\gamma^*$\cite{rugh1996linear}. Indeed, the output introduces a zero at the origin that
cancels a pole at the origin, rendering $x$ unobservable. 
Assuming no input $u$, we can write the system
(\ref{eq:simple-system}) as, \be{eq:simple-nonlin} \dot \xi = f(\xi),
\quad y = h(\xi), \ee where
$\xi = (x,\;z)^{\top},\;f = (z,\;-z)^{\top}$ and
$h(\xi) = \gamma(x)z$. We can construct the observation space,
$\mathcal{O}$ (set of all infinitesimal observables) by taking
$y = \gamma(x)z$ with all repeated time derivatives
\[
y^{(k)} = L_f^{(k)}(\gamma(x)z)
\]
as in \cite{mct,nijmeijer1990nonlinear}.  The
superscript ``$(k)$'' indicates $k$th order derivative. Note that
$L_f^{(k)}((\gamma(x)z))$ lies in the span of the functions
$\gamma^{(j)}(x)z^{j+1},\,j = 0,1,\ldots, k$. The rank condition on the
observability co-distribution \cite{mct,nijmeijer1990nonlinear} implies a sufficient condition for local
observability as follows
\cite{kunapareddy2018recovering}:
\be{}
    z^2(2(\gamma^{\prime}(x))^2-\gamma(x)\gamma^{\prime\prime}(x))\neq 0.
\ee 
Clearly for an non-hyperbolic $\gamma~(\neq1/(\alpha x+\beta)$, with
constants of integration $\alpha,\beta$), the non-zero velocity
requirement ($z\neq 0$) implies the need for active sensing to maintain the local
observability of the system~\cite{kunapareddy2018recovering,kunapareddyrecovering2016}.

What does this simplified example say about the need for active
sensing? As the proposition below illustrates, for such a system,
dynamic output feedback cannot asymptotically stabilize the origin
$(0,0)$, indicating the need for extra inputs to perform active state
estimation:

\bp{prop:nooutputfeedback} Consider the system
(\ref{eq:simple-system}). Let
\begin{equation}
  \label{eq:outputfeedback}
  \begin{aligned}
  \dot q &= g(q,y)\\
  u&=k(y,q)
\end{aligned}
\end{equation}
be a dynamic output feedback (Fig 1). Suppose
$(x^*,z^*,q^*)=(0,0,q^*)$ is an equilibrium of the coupled
system. Then all points $(\xi^*,0,q^*)$, $\xi^*\in\mathbb{R}$, are
equilibiria.\ep
\bpr
Since $(0,0,q^*)$ is an equilibrium, we see from the second equation
in (\ref{eq:simple-system}) that $k(0,q^*)=0$.  That means that
$k(\gamma(\xi^*)\cdot 0,q^*)=0$, i.e.\ $(\xi^*,0,q^*)$ is also an
equilibrium, for all $\xi^*\in\mathbb{R}$.
\epr

\begin{figure}[htbp!]
    \centering
    \includegraphics[width=5 cm]{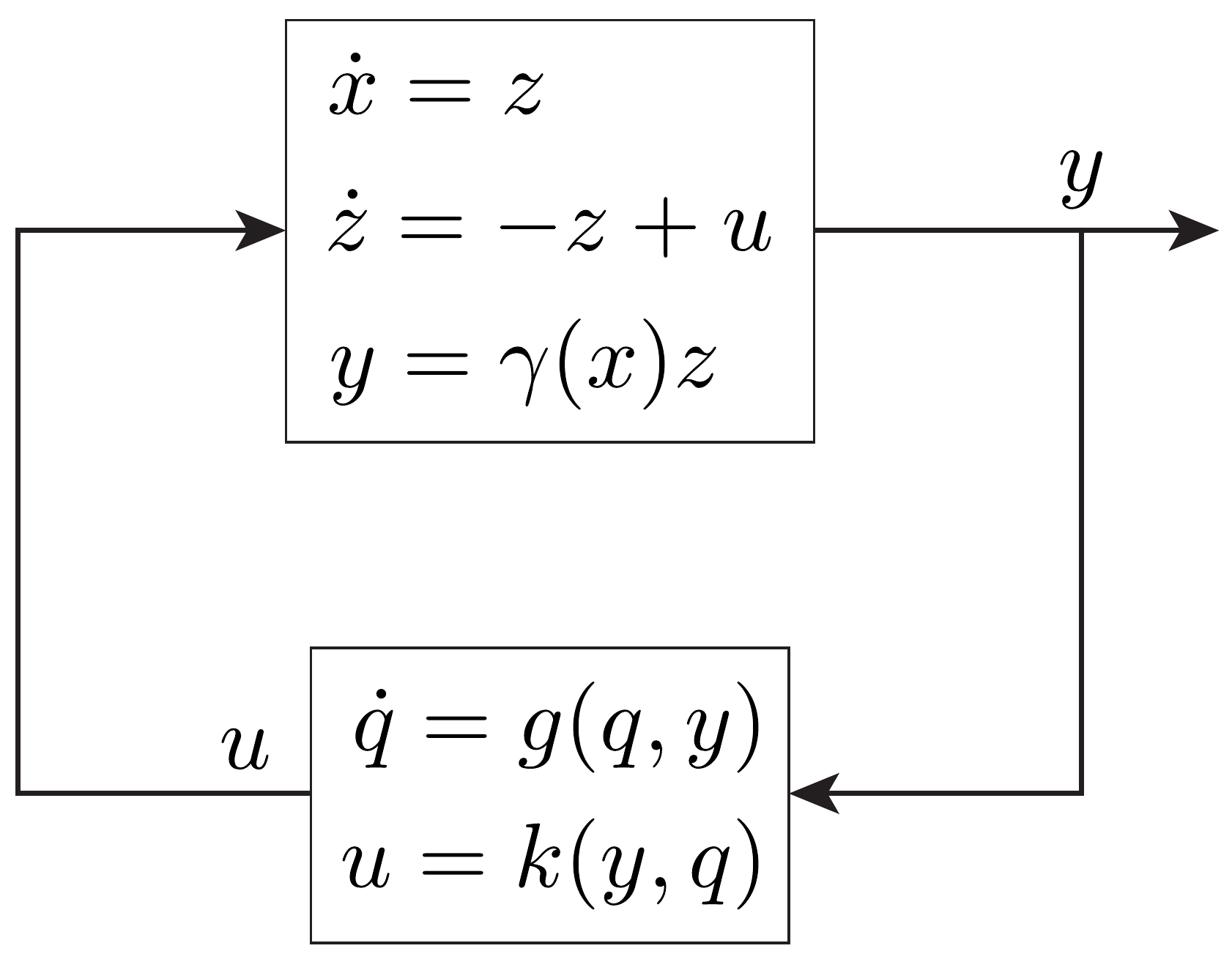}
    \caption{The system (1) cannot be stabilized to an equilibrium
      point by the dynamic feedback in (\ref{eq:outputfeedback}). }
    \label{fig:my_label}
\end{figure}

\br{remark:impossibility} The impossibility of stabilizing a system with
high-pass sensing to an equilibrium point, using only dynamic output
feedback, generalizes to the class of systems described below.
\er

\section{The class of systems and main result} \label{sec:mainresult}
Now we consider a general $2n$-dimensional, single-input, $n$-output system of the
following form: 
\beqn
\dot x &=& z\\
\dot z &=& F(z) \,+\, bu\\
y   &=& H(x)\,z
\eeqn
where the state space variable $\x\in \R^{2n}$ is partitioned as
$\x=(x^\top,z^\top)^\top$ into variables $x\in \R^n$ and $z\in \R^n$,
and $H(x)=\diag(\gamma _1(x_1), \ldots , \gamma _n(x_n))$ is a diagonal $n\times n$ matrix.
The entries of the column vector function $F(z)$ as well the functions
$\gamma _i(x_i)$ are real-analytic functions of their
arguments, and $b$ is a column vector of size $n$ all whose entries are nonzero.%

Let us introduce the following notations:
\beqn
f (\x)
&:=& \mypmatrix{z \cr F(z)}\\
g                    &:=& \mypmatrix{0 \cr b}\\
h (\x)
&:=& (h_1(\x),\ldots ,h_n(\x))^\top
\;=\;
(\gamma _1(x_1)\,z_1\ldots , \gamma _n(x_n)\,z_n)^\top\,,
\eeqn
so that, in standard form for nonlinear systems (see e.g.~\cite{mct}), our
system becomes
\be{eq:sys}
\dot \x = f(\x) + g u,\; y = h(\x) \,.
\ee
Given an input $u:[0,\infty )\rightarrow \R$ and an initial state $\x_0$, we denote by
$\varphi_{\x_0,u}(t)$ the solution of~\eqref{eq:sys}, that is,
$(d/dt)\varphi_{\x_0,u}(t) = f(\varphi_{\x_0,u}(t)) + g u(t)$ and $\varphi(0)=\x_0$.
Note that $\varphi_{\x_0,u}(t)$ is defined on some nontrivial time interval $[0,T)$
containing $t=0$.

\br{rem:technical}
Our results will hold in total generality, for inputs $u(\cdot )$ assumed to be
Lebesgue measurable locally bounded functions, in which case a solution is
an absolutely continuous function and the equality $\dot \x(t)=f(\x(t))+gu(t)$
holds almost-everywhere. If the reader prefers, one can restrict to inputs
$u$ which are piecewise continuous functions of time (with well-defined
one-sided limits at points of discontinuity) and solutions are piecewise
differentiable functions.
\er

We recall that two states $\x_0$ and $\hat\x_0$ are said to be
\textit{distinguishable} (by input/output measurements) if
there is some input $u(\cdot )$ and some time $t$ such that
$h(\varphi_{\x_0,u}(t))\not= h(\varphi_{\hat\x_0,u}(t))$, and that the system~\eqref{eq:sys}
is said to be \textit{observable} provided that every pair
of distinct states is distinguishable.
(An appendix reviews characterizations of observability in terms of Lie
derivatives, and these results are used in proofs.)

We will say that a function $\theta :\R\rightarrow \R$ is \textit{periodic} if there is
some nonzero $T\in \R$ (a ``period'')
such that $\theta (x)=\theta (x+T)$ for all $x\in \R$,
An \textit{aperiodic} function is one that is not periodic.

We will say that the system~\eqref{eq:sys} is \emph{aperiodic} if
none of the functions $\gamma _i$ are periodic.

Our main result is as follows:
\bt{theo:main}
The system \eqref{eq:sys} is observable if and only if it is aperiodic.
\et

The proof of Theorem follows immediately from Lemmas~\ref{lemma:aperiodic_observable} and~\ref{lemma:periodic_unobservable}.

\section{Proof of the main result} \label{sec:proof}

The key property that we need is as follows.

\bl{lemma:key_fact_lflgs}
For each integer $k\geq 0$, and each $i=1,\ldots ,n$, the following formulas hold:
\be{eq:lflg}
(L_fL_g)^kh_i (\x)\;=\; \gamma _i^{(k)}(x_i)\,b_i^k \,z_i 
\ee
\be{eq:lglflg}
L_g(L_fL_g)^kh_i (\x) \;=\; \gamma _i^{(k)}(x_i)\,b_i^{k+1} 
\ee
(where the superscript ``$(k)$'' indicates $k$th order derivative).
\el

\bpr
Fix any $i\in \{1,\ldots ,n\}$.
For $k=0$, formula~\eqref{eq:lflg} states $h_i (\x)=\gamma _i(x_i)z_i$.
which is true by definition of $h_i$.
Using induction, we will prove that, if formula~\eqref{eq:lflg} holds for a
given $k$ then~\eqref{eq:lglflg} holds for the same $k$, and~\eqref{eq:lflg}
holds for $k+1$.
Suppose that~\eqref{eq:lflg} is true.
Since
\[
L_g(L_fL_g)^kh_i \;=\;
\nabla [(L_fL_g)^kh_i] \cdot  g
\]
it follows that
\[
L_g(L_fL_g)^kh_i \;=\;
\nabla [\gamma _i^{(k)}(x_i)\,b_i^k \,z_i]\cdot  g \;=\;
\mypmatrix{\gamma _i^{(k+1)}(x_i)\,b_i^k\,z_i\,e_i\cr\gamma _i^{(k)}(x_i)\,b_i^k \,e_i}^\top \mypmatrix{0\cr b} \;=\; \gamma _i^{(k)}(x_i)\,b_i^{k+1}
\]
where $e_i$ is the $n$-vector with a ``1'' in position $i$ and zeroes elsewhere
.
Since
\[
(L_fL_g)^{k+1}h_i \;=\;L_f (L_g(L_fL_g)^kh_i)\;=\;\nabla [(L_g(L_fL_g)^kh_i] \cdot  f\,,
\]
we have that
\[
(L_fL_g)^{k+1}h_i \;=\;
\nabla [\gamma _i^{(k)}(x_i)\,b_i^{k+1}]\cdot  f \;=\;
\mypmatrix{\gamma _i^{(k+1)}(x_i)\,b_i^{k+1}e_i\cr 0}^\top \mypmatrix{z\cr F(z)}\;=\;
\gamma _i^{(k+1)}(x_i)\,b_i^{k+1}z_i\,.
\]
This completes the induction step.
\epr

\bl{lemma:periodic}
Suppose that $\gamma :\R\rightarrow \R$ is an aperiodic analytic function. Then, for each two
distinct $\ab,\bb\in \R$ there is some nonnegative integer $k$ such that
$\gamma ^{(k)}(\ab)\not= \gamma ^{(k)}(\bb)$.
\el

\bpr
Suppose, by way of contradiction, that there would be some pair $\ab\not= \bb$ for which
$\gamma ^{(k)}(\ab)=\gamma ^{(k)}(\bb)$ for all $k\geq 0$.
Let $T:=\bb-\ab$, so that $\bb=\ab+T$.
Consider the function $\beta (x):=\gamma (x+T)-\gamma (x)$, so that
\[
\beta ^{(k)}(\ab) \;=\; \gamma ^{(k)}(\ab+T)-\gamma ^{(k)}(\ab) \;=\; \gamma ^{(k)}(\bb)-\gamma ^{(k)}(\ab)
\;=\; 0
\]
for each $k\geq 0$.
Since the function $\beta $ is analytic, it follows that $\beta (x)\equiv 0$, which means
that $\gamma $ would be periodic of period $T$.
\epr

\bl{lemma:aperiodic_observable}
If the system~\eqref{eq:sys} is aperiodic,
then it is observable.
\el

\bpr
Consider two different states $\x_0$ and $\hat\x_0$. We need to show that
these two states are distinguishable.
We write these states in partitioned form as follows:
\[
\x_0 \;=\; \mypmatrix{x \cr z} \;=\; \mypmatrix{x_1\cr\vdots\cr x_n\cr z_1\cr\vdots\cr z_n}
\,,\quad
\hat\x_0 \;=\; \mypmatrix{\hat x \cr \hat z} \;=\; \mypmatrix{\hat x_1\cr\vdots\cr \hat x_n\cr \hat z_1\cr\vdots\cr \hat z_n}
\]
and consider two possible cases:
(a) $x=\hat x$, $z\not= \hat z$ and 
(b) $x\not= \hat x$.

Consider first case (a), so that $x_i=\hat x_i$ for all $i$. The set
\[
{\cal I}:=\{i\in \{1,\ldots ,n\} \st z_i\not= \hat z_i\}
\]
is nonempty.
Pick any $i\in {\cal I}$.
Since $\gamma _i$ is not periodic, it follows that 
$\gamma _i^{(k)}(x_i)\not= 0$ for some nonnegative integer $k$ (which may depend on $i$).
(If all derivatives were zero at a point, analyticity would imply $\gamma _i\equiv 0$, but the zero function is periodic.) 
For any such $i$ and $k$, 
$\gamma _i^{(k)}(x_i)=\gamma _i^{(k)}(\hat x_i)$ is nonzero and $z_i\not= \hat z_i$, so
$\gamma _i^{(k)}(x_i)\,b_i^kz_i\not=  \gamma _i^{(k)}(\hat x_i)\,b_i^k\hat z_i$.
From Equation~\eqref{eq:lflg} we have then that
\[
(L_fL_g)^kh_i(\x_0) 
\;=\; \gamma _i^{(k)}(x_i)\,b_i^kz_i
\;\not= \; \gamma _i^{(k)}(\hat x_i)\,b_i^k\hat z_i
\;=\; (L_fL_g)^kh_i(\hat\x_0)
\]
This means that the infinitesimal observable $(L_fL_g)^kh_i$ separates
the states $\x_0$ and $\hat\x_0$, so by Lemma~\ref{suft_elem_obs} these
states are distinguishable, as we claimed.

Next consider case (b).
Pick any $i\in \{1,\ldots ,n\}$ for which $x_i\not= \hat x_i$.
By Lemma~\ref{lemma:periodic}, there is some nonnegative integer $k$ such that
$\gamma ^{(k)}(x_i)\not= \gamma ^{(k)}(\hat x_i)$.
Then
\[
L_g(L_fL_g)^kh_i(\x_0) 
\;=\; \gamma _i^{(k)}(x_i)\,b_i^{k+1}
\;\not= \; \gamma _i^{(k)}(\hat x_i)\,b_i^{k+1}  
\;=\; L_g(L_fL_g)^kh_i(\hat\x_0)
\]
implies that the infinitesimal observable $L_g(L_fL_g)^kh_i$ separates
the states $\x_0$ and $\hat\x_0$, so by Lemma~\ref{suft_elem_obs} these
states are distinguishable, as we claimed.
\epr

\bl{lemma:periodic_unobservable}
If the system~\eqref{eq:sys} is not aperiodic,
then it is not observable.
\el

\bpr
Since~\eqref{eq:sys} is not aperiodic, there is some $i_0\in \{1,\ldots ,n\}$ and
some $T_{i_0}\not= 0$ such that $\gamma _{i_0}(x+T_{i_0})=\gamma _{i_0}(x)$ for all $x\in \R$.
Let $T$ be the $n$-vector that has this number $T_{i_0}$ as its $i_0$th
coordinate and zero in all other coordinates. Thus, $T\not= 0$, and
$\gamma _i(x+T_i)=\gamma _i(x)$ for all $i$ and all $x\in \R$.
Consider the following two distinct states $\x_0$ and $\hat\x_0$:
\[
\x_0 \;=\; \mypmatrix{0 \cr 0}\,,\quad
\hat\x_0 \;=\; \mypmatrix{T \cr 0}\,,\quad
\]
Consider any input $u(\cdot )$ and the respective solutions
$\x(t)=\varphi_{\x_0,u}(t)$ and $\hat\x(t)=\varphi_{\hat\x_0,u}(t)$.
In terms of the $x$ and $z$ components, we have that $z(t)=\hat z(t)$ for all
$t\geq 0$, because the $z$ component of the system does not depend on the $x$
component and the two initial conditions coincide on their $z$ components.
Thus, from
\[
x(t) \;=\; \int_0^\infty  z(t)\,dt
\]
and
\[
\hat x(t) \;=\; T + \int_0^\infty  \hat z(t)\,dt
\;=\; T + \int_0^\infty  z(t)\,dt
\]
we conclude that
\[
\hat x(t) \;=\; T + x(t)
\]
for all $t\geq 0$.
Since
$\gamma _i(x_i+T_i)=\gamma _i(x_i)$ for all $i$ and all $x_i\in \R$,
substituting $x_i=x_i(t)$ we have that
the output coordinates satisfy
\[
\hat y_i(t) \;=\; \gamma _i(\hat x_i(t)) \,\hat z_i(t)
\;=\;  \gamma _i(x_i(t)) \,z_i(t)
\;=\;  \gamma _i(x_i(t))\,z_i(t)
\;=\; y(t)
\]
for all $t\geq 0$.
This means that $\x_0$ and $\hat\x_0$ are indistinguishable.
\epr

\br{rem:no_analyticity}
We stated Lemmas~\ref{lemma:aperiodic_observable}
and~\ref{lemma:periodic_unobservable} separately because the latter one
does not require real-analyticity of any of the functions $\gamma _i$, and also
holds for a more general class of systems, namely
$\dot x = K(z)$, $\dot z = F(z,u)$.
No regularity properties whatsoever are required, except for
existence and uniqueness of solutions of the differential equations.
\er

\section{Conclusion}
Design of output feedback controllers commonly relies on the
separation principle, which allows designers to independently design
observers (based on sensor inputs) and controllers (designed assuming
full state measurements). In biological systems, this requirement for
separability may be violated. Specifically, high-pass sensing (which
we use to model adapting or perceptual fading sensory systems) causes
loss of observability for a class of systems. This manuscript presents
necessary and sufficient conditions for global observability for a
class of nonlinear systems with high-pass sensors. Though the system
structure was motivated by the locomotion dynamics of weakly electric
fish, it can be adopted to model behaviors of other animals with
translationally invariant plant dynamics and appropriately modeled
output measurements.

The goal of this work is to elucidate conditions that guarantee the
existence of inputs so that any two states can be distinguished, but a
characterization of ``good'' inputs that make the system
``sufficiently'' observable remains an open question.  Observability
and its dual, controllability are generic properties of a system,
although in practice are not always realizable due to practical
reasons, such as numerical conditioning \cite{cowan2012nodal}. Thus,
it remains unclear how best to design output feedback systems that
achieve a task-level control objective given the inseparability of
control and state estimation; a significant challenge in control
engingeering is to design a common framework to address both active
sensing and task-level control in a single design framework.

\appendix

\section{Relevant concepts related to nonlinear observability }

We review here a test for observability based on the  ``infinitesimal observables'' associated to a system.

Given any differentiable function $\alpha :\R^n\rightarrow \R$ and any vector field $X$, one
defines the \textit{Lie derivative of $\alpha $ along $X$} as the new function
with values
\[
L_g\alpha  \,:\, \R^n \rightarrow  \R \,:\, x \mapsto  \nabla \alpha (x) \cdot  X(x)
\]
where $\nabla\alpha $ is the gradient of $\alpha $ and ``$\cdot $'' indicates the dot or
inner product.
(This is the same as what in elementary calculus is called the
``directional derivative'' of the function $\alpha $ in the direction of $X$.)
This operation is multilinear: $L_{X+Y}\alpha =L_X\alpha +L_Y\alpha $ and
$L_{X}(\alpha +\beta )=L_X\alpha +L_x\beta $.

Consider a (generally multi-input multi-output) system in which inputs appear
linearly: 
\be{gen_nonlin_sys}
\dot  x \; = \; g_0(x) \,+\, \sum_{i=1}^m g_i(x) \,u_i
\ee
and with $p$ outputs $y_j=h_j(x)$, $j=1,\ldots ,p$
(in our application, $m=1$, $p=n$, $g_0=f$, $g_1=g$).
We only assume at first that all vector fields (that is, vector functions)
$g_i$ as well as the functions $h_j$ are infinitely differentiable; later
we impose real-analyticity (convergent power series around each state).

For any vector of nonnegative integers (though of as indices of vector fields)
\[
\mu  = (\mu _k, \ldots , \mu _1)\in \{0,\ldots ,m\}^k
\]
with $k\geq 1$, we define the \emph{infinitesimal observable} function
\be{elem-observables}
L_{\mu }h : \R^n \rightarrow  \R
\ee
by the iteration $L_{g_{\mu _k}}(L_{g_{\mu _{k-1}}}(\ldots (L_{g_{\mu _1}}h)\ldots ))$,
which we also write as $L_{g_{\mu _k}}L_{g_{\mu _{k-1}}}\ldots L_{g_{\mu _1}}h$.
We use power notation in the obvious form;
for example, $L_{g_1}^0h$ is just $h$, and $L_{g_1}^2h$ is the same as
$L_{(1,1)}h$, that is, $L_{g_1}(L_{g_1}h)$.
Finally, we let ${\cal O}$ be the set of all infinitesimal observables.

We say that two states $\x$ and $\hat\x$ are \textit{separated by ${\cal O}$} if there
exists some $\alpha \in {\cal O}$ such that $\alpha (\x)\not= \alpha (\hat\x)$. 
We have the following well-known fact (see e.g.~\cite{mct} and references
there), which we prove here for ease of reference.

\bl{suft_elem_obs}
If two states are separated by ${\cal O}$, then they are distinguishable.
\el

\bpr
We prove the contrapositive: if two states $\x$ and $\hat\x$ are not
distinguishable, then they cannot be separated by ${\cal O}$, that is,
$\alpha (\x)=\alpha (\hat\x)$ for all $\alpha \in {\cal O}$. 
Pick $\x$ and $\hat\x$ that are not distinguishable,
and consider a piecewise constant input on $k$ intervals:
$u(\cdot )$ has a constant value $u^1$ on an interval $[0,t_1)$, a constant value
$u^2$ on $[t_1,t_1+t_2)$, \ldots , and $u^k$ on $[t_1+\ldots +t_{k-1},t_1+\ldots +t_k)$.
For small enough $t_i$'s there is a solution of the differential equation
from both initial conditions $\x$ and $\hat\x$.
Since these two states are indistinguishable, the resulting output at time
$t=t_1+\ldots +t_k$ is the same, when starting from either initial state.
In general, let us denote the $j$th coordinate of this output value by
\be{output-pc}
h_j(t_1,t_2,\ldots ,t_k,u^1,u^2,\ldots ,u^k,\widetilde \x)
\ee
when the initial state is $\widetilde \x$.
It follows that the derivatives with respect to the $t_i$'s of this output are
also equal, for $\x$ and $\hat\x$, for every such piecewise constant input.
One may prove by induction that
\[
\frac{\partial ^k}{\partial t_1 \ldots  \partial t_k}
\bigg| _{t_1=t_2=\ldots =0}
h_j  (t_1,t_2,\ldots ,t_k,u^1,u^2,\ldots ,u^k,\widetilde \x) 
\;=\;
L_{X_1} L_{X_2} \ldots  L_{X_k} h_j(\widetilde \x)
\]
where $X_l(x) = g_0(x) + \sum_{i=1}^m u_i^l g_i(x)$.
In summary,
\[
L_{X_1} L_{X_2} \ldots  L_{X_k} h_j(\x)
\;=\;
L_{X_1} L_{X_2} \ldots  L_{X_k} h_j(\hat\x)
\]
for any $k$ and any vector $(u^1,\ldots ,u^k)\in \R^k$.
Using multilinearity, $L_{X_1} L_{X_2} \ldots  L_{X_k} h_j(x)$ can be expanded as a
polynomial on the $u^1,\ldots ,u^k$ whose coefficients are exactly the elementary
observables.
For example, if $k=2$ and $m=1$,
\beqn
L_{X_1} L_{X_2} h_j &=&
L_{g_0+u^1 g_1}(L_{g_0+u^2 g_1} h_j)\\
&=&
L_{g_0}(L_{g_0+u^2 g_1} h_j)
+
u^1L_{g_1(x)}(L_{g_0+u^2 g_1} h_j)\\
&=&
L_{g_0}(L_{g_0}h_j + u^2 L_{g_1} h_j)
+
u^1L_{g_1(x)}(L_{g_0}h_j + u^2 L_{g_1}h_j)\\
&=&
L_{g_0}^2h_j +
u^1 L_{g_1}L_{g_0}h_j + 
u^2 L_{g_0}L_{g_1}h_j +
u^1u^2 L_{g_1}^2h_j\,.
\eeqn
Since two polynomial functions are equal if and only if their coefficients of
equal powers are equal, it follows that $L_{\mu }h_j(\x)=L_{\mu }h_j(\hat\x)$ for
all $k$ and all indices $\mu $, and hence these states cannot be separated by
${\cal O}$.
\epr

We'll say that \textit{observables separate states} if any two states can be
separated by ${\cal O}$. A consequence of the above is:

\bc{cor:observable}
If observables separate states, then the system is observable.
\ec

To rephrase the corollary, a sufficient condition for observability is that
the mapping 
\[
\x \;\mapsto  \; \{\alpha (\x), \alpha \in {\cal O}\} \,,
\]
which sends each state into the infinite sequence of possible observables
evaluated at that state, be one-to-one.

This condition is also necessary when all the functions are real-analytic:

\bt{theo:real-analytic-observables}
Suppose that the vector fields $g_i$ as well as the functions $h_j$ are real
analytic. Then, the following two properties are equivalent:
\ben
\item
  The system is observable.
\item
  Observables separate states.
\een
\et

\bpr
One implication is given by Corollary~\ref{cor:observable}.
To prove the converse, we need to show that if observables do not separate
states then the system is not observable.
Indeed, suppose that there is some pair of states $\x$ and $\hat\x$ such that
$\alpha (\x)=\alpha (\hat\x)$ for all $\alpha \in {\cal O}$. We want to show that these two states are
not distinguishable.
We first show that these two states are not distinguishable by means of piecewise
constant inputs.
This follows from the construction in Lemma~\ref{suft_elem_obs}. For any
given piecewise constant input, we have that
\beqn
\frac{\partial ^k}{\partial t_1 \ldots  \partial t_k}
\bigg| _{t_1=t_2=\ldots =0}
h_j  (t_1,t_2,\ldots ,t_k,u^1,u^2,\ldots ,u^k,\x)
&=&
L_{X_1} L_{X_2} \ldots  L_{X_k} h_j(\x)\\
\frac{\partial ^k}{\partial t_1 \ldots  \partial t_k}
\bigg| _{t_1=t_2=\ldots =0}
h_j  (t_1,t_2,\ldots ,t_k,u^1,u^2,\ldots ,u^k,\hat\x) &=&
L_{X_1} L_{X_2} \ldots  L_{X_k} h_j(\hat\x)
\eeqn
and the two right-hand sides coincide because observables do not separate these two states.
Now, the maps
\[
h_j  (t_1,t_2,\ldots ,t_k,u^1,u^2,\ldots ,u^k,\widetilde \x)
\]
are analytic functions of the times $t_i$ (see e.g.~\cite{mct}), and two
analytic functions that have the same derivatives at one point must be the
same, so 
\[
h_j  (t_1,t_2,\ldots ,t_k,u^1,u^2,\ldots ,u^k,\x) \;=\;
h_j  (t_1,t_2,\ldots ,t_k,u^1,u^2,\ldots ,u^k,\hat\x)
\]
for any such piecewise constant input.
To finalize the proof, observe that piecewise constant inputs are dense in
the set of measurable inputs (see e.g.~\cite{mct}, Remark C.1.2), and that the
state (and hence output) is continuous with respect to the weak topology on
inputs (see e.g.~\cite{mct}, Theorem 1).
Thus the states $\x$ and $\hat\x$ are not distinguishable.
\epr

\br{rem:real-analytic-observables}
There is a completely different proof of the same fact, using Corollary 5.1 in
\cite{MR1025981}. This Corollary says that two states are indistinguishable if
and only if for every polynomial-in-time input, the derivatives at time zero
of the output $y(t)$ are the same. Since polynomial inputs are dense in the
set of all inputs, the result follows by a density argument.
\er

\section{Relevant concepts related to local nonlinear observability }

We review here a test for local observability based on the observation space, $\mathcal{O}$.
For a general autonomous multi-output system (with $u=0$ in eq. \ref{gen_nonlin_sys})
\begin{equation}
\begin{aligned}
\label{gen_autonom_sys}
\dot  x \; &= \; g_0(x),\quad x\in X \\
y_j \; &= h_j(x),\quad j\in p\;
\end{aligned}
\end{equation}
one can construct $\mathcal{O}$ using $p$ outputs, $y_j=h_j(x)$ and the repeated time derivatives, $y_j^{(k)}=L_{g_0}^{(k)}h_j(x)$. Consequently, one can define the observability codistribution, $d\mathcal{O}$ as
\[
d\mathcal{O}(\x) = \spn\bigg\{\dfrac{\partial \alpha}{\partial x} \,|\,\alpha\in \mathcal{O}\bigg\}, \quad \x \in X,
\]

\bt{theo:local-observ-rank-cond}
The system (\ref{gen_autonom_sys}) is locally observable at $\x_0$ if 
\beqn
\dim d\mathcal{O}(\x_0) = \dim (X).
\eeqn
\et
The proof of the theorem \ref{theo:local-observ-rank-cond} is given in
\cite{mct} (Remark 6.4.2, pp. 281-282) or
\cite{nijmeijer1990nonlinear} (pp. 95-96).

\bibliographystyle{ieeetr}
\bibliography{observability_cascade}

\end{document}